\documentclass[preprint,review,12pt]{elsarticle}

\usepackage{amssymb}

\journal{Journal}
\begin{document}

\begin{frontmatter}
\title{On the Firoozbakht's conjecture}
\author{Ahmad Sabihi\fnref{fn1}} 
\ead{sabihi2000@yahoo.com}
\address{Teaching professor and researcher at some universities of Iran}
\fntext[fn1]{ Fourth Floor, Building 30, Poormehrabi Alley, Mofatteh St., Bozorgmehr Ave., Isfahan, Iran.} 
\begin{abstract}
This paper proves Firoozbakht's conjecture using Rosser and Schoenfelds' inequality on the distribution of primes. This inequality is valid for all natural numbers ${n\geq 21}$. Firoozbakht's conjecture states that if $ {p_{n}}$ and  ${p_{(n+1)}}$ are consecutive prime numbers, then ${p_{(n+1)}^{1/(n+1)}< p_{n}^{1/n}}$ for every ${n\geq 1}$. Rosser's inequality for the ${n}$th and ${(n+1)}$th roots, changes from strictly increasing to strictly decreasing for ${n\geq 21}$. The inequality is considered for ${n>e^{e^{3/2}}}$, i.e., ${n\geq 89}$, but since the inequalities for ${n\geq 195340>e^{e^{5/2}}}$, are also required, these inequalities are explicitly proven as well. Silva has already verified Firoozbakht's conjecture up to ${p_{n}<4 \times 10^{18}}$, and the additional theorem is proven here that there is the smallest natural number, ${m>n\geq 1}$ and ${p_{m}^{1/m}< p_{n}^{1/n}}$. It is also shown that there is a unique one to one function, which maps each element ${p_n}$ to each element ${p_{n}^{1/n}}$ for every ${n\geq 1}$ and ${1<p_{n}^{1/n}\leq 2}$. Finally, it is proved that there is a strictly decreasing sequence and Firoozbakht's conjecture would be true for all ${n\geq 1}$.
\end{abstract}

\end{frontmatter}

\section{Introduction}

Firoozbakht's conjecture was proposed by the mathematician Farideh Firoozbakht in 1982 \cite{Fe}, and relates to the distribution of primes, specifically Cramer primes. Cramer's conjecture \cite{C} states that gaps between consecutive prime numbers can have a supermum 1 with regard to ${(\log p_{n})^{2}}$ ($\log$ refers to natural logarithm throughout the paper) as ${\lim_{n \rightarrow \infty}sup\frac{(p_{n+1}-p_{n})}{(\log p_{n})^{2}}=1}$. Shanks' conjecture \cite{S} (( ${p_{n+1}-p_{n})\sim (\log p_{n})^{2}}$) gives a somewhat stronger statement than Cramer's, and Nicely \cite{N} provided many calculations on the prime gaps and their relationship to Cramer's conjecture. Firoozbakht's conjecture has been verified up to ${4.444 \times 10^{12}}$ \cite{F}, and ${4 \times 10^{18}}$ \cite{T}. Kourbatov \cite{K} showed that if Firoozbakht's conjecture is true for the ${k}$th prime, then the inequality ${p_{k+1}-p_{k}<(\log p_{k})^{2}-\log p_{k}-b}$ is also true for ${b=1}$ and with ${b=1.17}$ implies Firoozbakht's conjecture. Rosser and Schoenfeld \cite{RS} proved the inequalities \begin{itemize}\item For ${x\geq 59}$ 
\begin{equation}
{\frac{x}{\log x}(1+\frac{1}{2\log x})<\pi(x)<\frac{x}{\log x}(1+\frac{3}{2\log x})}.
\end{equation}

\item For ${x\geq 41}$ 
\begin{equation}
{x(1-\frac{1}{\log x})<\theta(x)<x(1+\frac{1}{2\log x})}.
\end{equation}

\item For ${x\geq 121}$ 
\begin{equation}
{x(1-\frac{1}{\log x}+\frac{0.98}{x^{0.5}})<\psi(x)<x(1+\frac{1}{2\log x}+\frac{1.02}{x^{0.5}}+\frac{3}{x^{2/3}})}.
\end{equation}

\item For all the natural numbers ${n\geq 21}$
\begin{equation}\label{eq1.4}
{n(\log n+\log \log n-3/2)<p_{n}<n(\log n+\log \log n-1/2)}.
\end{equation}
\end{itemize}

We make use of (\ref{eq1.4}) throughout the paper to prove Firoozbakht's conjecture. In Section~\ref{sc2} we prove Theorem~1 and three related lemmas to show that there exist real numbers ${p_{m}^{1/m}}$ which can be of values less than ${p_{n}^{1/n}}$ for ${n\geq 1}$. We subsequently develop three inequalities based on (1.4). In Section~\ref{sc3} we prove the inequalities of Section~\ref{sc2}, to finalize the proof of Theorem~1. In Section~\ref{sc4}, we prove Theorem~2 on a sequence of real numbers, which implies Firoozbakht's conjecture.

\setcounter{equation}{0}
\section{Methods and Statement of the Purpose} \label{sc2}

\subsection{Theorem 1}

\textit{Let ${n\geq 1}$ be a natural number. There is the smallest natural number, ${m}$, such that ${m>n\geq 1}$ and ${p_{m}^{1/m}< p_{n}^{1/n}}$, where ${p_{n}}$ and ${p_{m}}$ are the ${n}$th and ${m}$th prime numbers, respectively. Conversely, ${p_{m}^{1/m}< p_{n}^{1/n}}$ implies $m>n$}.

\textbf{\textit{Proof}}

Rosser and Schoenfelds' inequality on the distribution of primes is as follows:

If ${p_{n}}$ is the ${n}$th prime number, then for ${n\geq 21}$,
\begin{equation}
{n(\log n+\log \log n-3/2)<p_{n}<n(\log n+\log \log n-1/2)},
\end{equation}
and the $n$th root of the inequality is
\begin{equation}\label{eq2.2}
{\{n(\log n+\log \log n-3/2)\}^{1/n}<p_{n}^{1/n}<\{n(\log n+\log \log n-1/2)\}^{1/n}}.
\end{equation} 

A similar inequality can be established for the ${(n+1)}$th prime number,
\begin{eqnarray}
{\{(n+1)(\log (n+1)+\log \log (n+1)-3/2)\}^{1/(n+1)}<p_{(n+1)}^{1/(n+1)}}<\nonumber\\
 {\{(n+1)(\log (n+1)+\log \log (n+1)-1/2)\}^{1/(n+1)}}
\label{eq2.3}
\end{eqnarray} 

We want to prove that left and right sides of inequality (\ref{eq2.3}) are less than the corresponding sides of inequality (\ref{eq2.2}). Comparing the right side of (\ref{eq2.3}) with (\ref{eq2.2}), we find
\begin{eqnarray}
{\{(n+1)\log (n+1)+(n+1)\log \log (n+1)-(n+1)/2)\}^{1/(n+1)}}< \nonumber\\ 
{\{n\log n+n\log \log n-n/2)\}^{1/n}}, \label{eq2.4}
\end{eqnarray}
for ${n\geq 89}$, and the corresponding left sides is 
\begin{eqnarray}
{\{(n+1)\log (n+1)+(n+1)\log \log (n+1)-3(n+1)/2)\}^{1/(n+1)}}<\nonumber\\ 
{\{n\log n+n\log \log n-3n/2)\}^{1/n}}, \label{eq2.5}
\end{eqnarray}
for ${n\geq 195340}$. Thus, 
\begin{eqnarray}
{p_{m}^{1/m}<\{m\log m+m\log \log m-m/2\}^{1/m}<}\nonumber\\ {\{n\log n+n\log \log n-3n/2\}^{1/n}< p_{n}^{1/n}}  \label{eq2.6}
\end{eqnarray}
for $m\geq 3n$ and  $n\geq 195340$.

If (\ref{eq2.4}) can be proven for ${n\geq 89}$, then (\ref{eq2.5}) is also proven for ${n\geq 195340}$. Therefore, it is sufficient that we only prove (\ref{eq2.4}) and (\ref{eq2.6}). Note that (\ref{eq2.6}) shows that inequality ${p_{m}^{1/m}< p_{n}^{1/n}}$ holds for all  $n\geq 195340$ .

We present several lemmas, which along with the inequalities (2.4), (2.5) and (2.6)  prove Theorem~1.

\subsubsection{Lemma 1}

\textit{If ${n\geq 21}$, then ${(n+1)\log(n+1)<n^{(1+1/n)} (\log n)^{(1+1/n)}}$ holds.}

\textbf{\textit{Proof}}

The proof consists of two steps:
\begin{enumerate}
\item ${(n+1)<n^{(1+1/n)}}$

Note that ${(1+\frac{1}{n})^{n}<e<3}$ holds for all ${n\geq 1}$ because the sequence ${(1+\frac{1}{n})^{n}}$ acts as a strictly increasing sequence for ${n\geq 1}$. Therefore, $${(1+\frac{1}{n})^{n}<3\leq n}$$
for ${n\geq 3}$.
Taking the $n$th root, ${(1+\frac{1}{n})<n^{1/n}}$, then ${(n+1)<n^{(1+1/n)}}$ 
for ${n\geq 21>3}$.

\item ${\log(n+1)< (\log n)^{(1+1/n)}}$

Since ${(1+\frac{1}{n})^{n}\leq e<3}$, if we assume ${\log n>3}$ or ${n>e^{3}}$, then ${(1+\frac{1}{n})^{n}<3}{<\log n}$ and ${(1+\frac{1}{n})^{n}<\log n}$. Taking the ${n}$th root, we have  ${(1+\frac{1}{n})<}$${(\log n)^{1/n}}$ and ${n+1<n(\log n)^{1/n}}$, and so ${1<n(\log n)^{1/n}-n}$. Multiplying both sides of the inequality by ${\log n}$, we have
$${1<\log n<n\log n((\log n)^{1/n}-1)}$$ 
for ${n>e^{3}}$, and taking the exponent, $${(1+\frac{1}{n})^{n}}<{e}<{e^{\log n}< e^{n\log n((\log n)^{1/n}-1)}}.$$
The ${n}$th root shows that $${(1+\frac{1}{n})<}{e^{\log n((\log n)^{1/n}-1)}},$$ 
and taking the logarithm, we have
$${\log n+\log(1+1/n)}{<(\log n)^{(1+1/n)}},$$ 
so that 
$${\log(n+1)<(\log n)^{(1+1/n)}}$$ 
for ${n\geq 21>e^{3}}$.
\end{enumerate}

Multiplying both sides of the step 1 and 2 inequalities proves Lemma~1.

\subsubsection{Lemma 2}

\textit{Let ${n\geq 1}$, then the inequality ${1<p_{n}^{1/n} \leq 2}$ holds for all natural numbers.}

\textbf{\textit{Proof.}}
\begin{enumerate}
\item We show that ${p_{n}^{1/n} \leq 2}$. Note that ${p_{n}\leq 2^{n}}$ is true for all natural numbers by induction (if ${n=1}$, ${p_{1}=2 \leq 2}$). If the inequality ${p_{n}\leq 2^{n}}$ is true for ${n}$,then we must prove that it is also true for ${p_{(n+1)}\leq 2^{(n+1)}}$. 

Combining Bertrand's postulate, ${p_{n}\leq p_{(n+1)}\leq 2p_{n}}$, and our assumption, ${p_{(n+1)}\leq 2p_{n}\leq 2^{(n+1)}}$, and so  ${p_{(n+1)}\leq 2^{(n+1)}}$. Thus ${p_{n}\leq 2^{n}}$ holds for all the natural numbers and implies that ${p_{n}^{1/n} \leq 2}$ also holds.

\item We show that ${1<p_{n}^{1/n}}$ holds for all natural numbers. Trivially, because of the monotony of $n$th radical, we have

$1<p_{n}$ and taking $n$th radical, $\sqrt[n]{1}<\sqrt[n]{p_{n}}$ , and so ${1<p_{n}^{1/n}}$ for $n\geq 1$
\end{enumerate}

\subsubsection{Lemma 3} 

\textit{Let ${f(n)=p_{n}^{1/n}}$ for ${n\geq 1}$, then there is a one to one function on the natural numbers set such that ${f(n):\Bbb N \mapsto ]1,2]}$, where $\Bbb N$ is the set of natural numbers}.

\textbf{\textit{Proof.}} 

In Lemma~2, we showed that ${p_{n}^{1/n}\in ]1,2]}$ for ${n\geq 1}$. Let ${n}$ and ${m}$ be natural numbers such that ${m\geq n \geq1}$, ${f(n)=p_{n}^{1/n}}$, and ${f(m)=p_{m}^{1/m}}$. If they contradict the property of a one to one function, then ${f(n)=p_{n}^{1/n}=f(m)=p_{m}^{1/m}}$ and ${n\neq m}$, which would mean that ${p_{n}=p_{m}^{n/m}}$. However, ${p_{n}}$ is a prime number, and, since this equation must hold, ${0<n/m\leq 1}$, and the only solution is ${n=m}$, which completes the proof. 

Thus ${f(n)}$ is a one to one function.

\setcounter{equation}{0}
\section{Inequalities (\ref{eq2.4}), (\ref{eq2.5}) and (\ref{eq2.6})} \label{sc3}

As stated in Section~\ref{sc2}, to prove Theorem~1 it is sufficient to prove the (\ref{eq2.4}) and (\ref{eq2.6}). Inequality (\ref{eq2.5}) is concluded from (\ref{eq2.4}) since ${n\geq 195340}$.

\subsection{Proof of inequality (\ref{eq2.4})}

If we take the right and left sides of (\ref{eq2.4}) to the ${n(n+1)}$th power, then
\begin{eqnarray}
(n+1)^{n}\{\log (n+1)+\log \log (n+1)-1/2\}^{n}< \nonumber\\ 
n^{(n+1)}\times \{\log n+\log \log n-1/2\}^{(n+1)} \label{eq3.1} 
\end{eqnarray}

The proof comprises two steps:
\begin{enumerate}
\item ${(n+1)^{n}<n^{(n+1)}}$ 

From Lemma~1, step 1, ${(n+1)<n^{(1+1/n)}}$ for all ${n\geq 3}$. Therefore, if one takes each side of the inequality to the ${n}$th power, then the result is as above.

\item ${\{\log (n+1)+\log \log (n+1)-1/2\}^{n}}<{\{\log n+\log \log n-1/2\}^{(n+1)}}$ 

To prove this inequality, we employ binomial expansion of each side and show that each of the terms of the expanded left side is less than the corresponding right side term. We then show that the ${(n+1)}$th term of the left side is less than the sum of the expanded ${(n+1)}$th and ${(n+2)}$th terms of the right of the inequality.

The ${i}$th term of the expanded left and right sides of the inequality of step 2 (${i=1,\dots, n}$) is the inequality 
\begin{eqnarray}
 \left(%
\begin{array}{cc}  
{n}\\ 
{ i-1}\\
\end{array}%
\right){(\log(n+1))^{n-(i-1)}\{\log \log (n+1)-\frac{1}{2}\}^{(i-1)}}< \nonumber\\ \left(%
\begin{array}{cc}   
{ n+1}\\ 
{ i-1}\\
\end{array}%
\right){(\log n)^{n+1-(i-1)}\{\log \log n-\frac{1}{2}\}^{(i-1)}}, \label{eq3.2}
\end{eqnarray}
where $\left(
\begin{array}{cc}   
 {n}\\ 
{i-1}\\
\end{array}%
\right)$ and $\left(%
\begin{array}{cc}   
 {n+1}\\ 
{i-1}\\
\end{array}%
\right)$ are permutations of ${(i-1)}$ of ${n}$ and ${n+1}$, respectively.
\end{enumerate}

Proof of (\ref{eq3.2}) consists of three steps
\begin{enumerate}
\item $\left(%
\begin{array}{cc}   
 {n}\\ 
{i-1}\\
\end{array}%
\right)\leq \left(%
\begin{array}{cc}   
{n+1}\\ 
{i-1}\\
\end{array}%
\right)$ for ${i=1,\dots,n}$. 

The proof can be done using Pascal's identity

$\left(%
\begin{array}{cc}   
 {n}\\ 
{i-2}\\
\end{array}%
\right)+\left(%
\begin{array}{cc}   
{n}\\ 
{i-1}\\
\end{array}%
\right)=\left(%
\begin{array}{cc}   
{n+1}\\ 
{i-1}\\
\end{array}%
\right)$ for ${i=1,\dots,n}$. 

Therefore, $\left( 
\begin{array}{cc}   
 {n}\\ 
{i-1}\\
\end{array}%
\right)\leq \left(%
\begin{array}{cc}   
 {n+1}\\ 
{i-1}\\
\end{array}%
\right)$, which completes the proof.

\item The second portion of (\ref{eq3.2})holds for ${i=1, \dots, n}$.

The second part of (\ref{eq3.2}) may be expressed as
\begin{eqnarray}
(\log(n+1))^{n-(i-1)}\{\log \log (n+1)-1/2\}^{(i-1)}< \nonumber\\ 
(\log n)^{n+1-(i-1)}\times \{\log \log n-1/2\}^{(i-1)}. \label{eq3.3}
\end{eqnarray}
For ${i=1}$, (\ref{eq3.3}) becomes ${(\log(n+1))^{n}<(\log n)^{n+1}}$. From Lemma~1, step~2, ${\log(n+1)< (\log n)^{(1+1/n)}}$ for ${n\geq 21}$. Therefore, taking each side to the ${n}$th power, we have the required result.

We now prove that (\ref{eq3.3}) is satisfied for ${i=2, \dots, n}$. 
Let ${\log (n+1)=\log n\{1+\frac{\log(1+1/n)}{\log n}\}}$, then (\ref{eq3.3}) becomes 
\begin{eqnarray}
{(\log n)^{(n-(i-1))}\{1+\frac{\log(1+1/n)}{\log n}\}^{(n-(i-1))}\{\log \log (n+1)-1/2\}^{(i-1)}}\nonumber\\< {(\log n)^{(n-(i-2))}\{\log \log n-1/2\}^{(i-1)}}. \label{eq3.4}
\end{eqnarray}
Eliminating the common powers of ${\log n}$ from both sides, 
\begin{eqnarray}
\{1+\frac{\log(1+1/n)}{\log n}\}^{(n-(i-1))}\{\log \log (n+1)-1/2\}^{(i-1)} \nonumber\\ 
<(\log n)\times \{\log \log n-1/2\}^{(i-1)}, \label{eq3.5}
\end{eqnarray}
and dividing both sides by ${(\log \log n-\frac{1}{2})^{(i-1)}}$, 
\begin{equation} \label{eq3.6}
{\{1+\frac{\log(1+1/n)}{\log n}\}^{(n-(i-1))}(\frac{\{\log \log (n+1)-1/2\}}{\{\log \log n-1/2\}})^{(i-1)}<\log n}
\end{equation}

The first term of the left side of (\ref{eq3.6}) is less than ${e^{\frac{1}{\log n}}}$.

Since ${(1+1/n)^{(n-(i-1))}=\frac{(1+1/n)^{n}}{(1+1/n)^{(i-1)}}}$, ${1<(1+1/n)^{(n-(i-1))}<e}$ and $${\log (1+1/n)^{(n-(i-1))}<\log e=1}$$,then, 
$${[1+\frac{\log(1+\frac{1}{n})}{\log n}]^{(n-(i-1))}= \{[1+\frac{\log(1+1/n)}{\log n}]^{\frac{\log n}{\log(1+1/n)}}\}^{\frac{\log (1+1/n)^{(n-(i-1))}}{\log n}}<e^{\frac{1}{\log n}}}$$ 

The second term of the left side of (\ref{eq3.6}) is less than ${e^{\frac{(i-1)\log(1+1/n)}{\log n(\log \log n-1/2)}}}$.

Let 
$${\{\log \log (n+1)-1/2\}=\log \log n+\log[1+\frac{\log(1+1/n)}{\log n}]-1/2},$$
 where $S = {\{1+\frac{\log[1+\frac{\log(1+1/n)}{\log n}]}{\{\log \log n-1/2\}}\}}$, and $T = {\{1+\frac{\log(1+1/n)}{\log n}\}}$. Then, 
$$
{S^{(i-1)}=(\frac{\{\log \log (n+1)-1/2\}}{\{\log \log n-1/2\}})^{(i-1)}=\{1+\frac{\log[1+\frac{\log(1+1/n)}{\log n}]}{\{\log \log n-1/2\}}\}^{(i-1)}}=$$  $${\{S^{\frac{\{\log \log n-1/2\}}{\log[1+\frac{\log(1+1/n)}{\log n}]}}\}^{\{(i-1)\frac{\log[1+\frac{\log(1+1/n)}{\log n}]}{\{\log \log n-1/2\}}\}}}
$$ 
It is trivial that if ${x>0}$, then ${\log (1+x)<x}$. Thus, ${\log[1+\frac{\log(1+1/n)}{\log n}]<\frac{\log(1+1/n)}{\log n}}$, and  ${S^{(i-1)}<e^{\frac{(i-1)\log(1+1/n)}{\log n(\log \log n-1/2)}}}$.
Thus, ${T^{(n-(i-1))}S^{(i-1)}<e^{\frac{1}{\log n}}e^{\frac{(i-1)\log(1+1/n)}{\log n(\log \log n-1/2)}}}$ for ${n\geq 89}$. 

Therefore, the sum of the powers is
$$
{\frac{(i-1)\log(1+1/n)}{(\log n)(\log \log n-1/2)}+\frac{1}{\log n}<\frac{(n-1)\log(1+1/n)}{(\log n)(\log \log n-1/2)}}+\\
\frac{1}{\log n}
$$
for $ i=2,\dots,n$. Since 
$${(n-1)\log(1+1/n)=\log(1+1/n)^{(n-1)}<\log(1+1/n)^{n}<1},$$ 

$${\frac{(n-1)\log(1+1/n)}{(\log n)(\log \log n-1/2)}+\frac{1}{\log n}<\frac{1}{(\log n)(\log \log n-1/2)}+\frac{1}{\log n}}.$$
However, for ${n\geq 89}$,
$$
{\frac{1}{(\log n)(\log \log n-1/2)}+\frac{1}{\log n}<\frac{1}{(\log 89)(\log \log 89-1/2)}+\frac{1}{\log 89}\approx 0.45}$$ $${<{\log \log 89}<{\log \log n}}.
$$
Thus, 
$$
{\frac{(i-1)\log(1+1/n)}{(\log n)(\log \log n-1/2)}+\frac{1}{\log n}<\log \log n},
$$
and 
$$
{T^{(n-(i-1))}S^{(i-1)}=\{1+\frac{\log(1+1/n)}{\log n}\}^{(n-(i-1))}(\frac{\{\log \log (n+1)-1/2\}}{\{\log \log n-1/2\}})^{(i-1)}<}$$  $${e^{(\frac{1}{\log n})}e^{\{\frac{(i-1)\log(1+1/n)}{\log n(\log \log n-1/2)}\}}<e^{\log \log n}=\log n}.$$

This completes the proof of (\ref{eq3.6}) and consequently (\ref{eq3.3}) for ${i=2, \dots, n}$ and ${n\geq 89}$.

\item We prove inequality (\ref{eq3.2}) holds when the ${i=(n+1)}$th term of left side of and the sum of the ${i=(n+1)}$th and ${i=(n+2)}$th terms of the right side are compared. 

This means that we need to show that 
\begin{eqnarray}
{\{\log \log(n+1)-1/2\}^{n}<(n+1)\log n(\log \log n-1/2)^{n}+}\nonumber\\ 
{(\log \log n-1/2)^{n+1}}. \label{eq3.7}
\end{eqnarray}
holds for ${n\geq 89}$.

Inequality (\ref{eq3.7}) can be rewritten as 
\begin{eqnarray}
{\{\log \log(n+1)-1/2\}^{n}<(\log \log n-1/2)^{n}\times} \nonumber\\ {\{(n+1)\log n+\log \log n-1/2\}}, \label{eq3.8}
\end{eqnarray}
and taking the ${n}$th root, 
\begin{eqnarray}
{\{\log \log(n+1)-1/2\}<(\log \log n-1/2)\times} \nonumber\\ {\{(n+1)\log n+\log \log n-1/2\}^{1/n}}. \label{eq3.9}
\end{eqnarray}

Dividing both sides by ${(\log \log n-1/2)}$, and taking the factor ${(n+1)^{1/n}(\log n)^{1/n}}$, 
$$
{\frac{(\log \log(n+1)-1/2)}{(\log \log n-1/2)}<(n+1)^{1/n}(\log n)^{1/n}\{1+\frac{(\log \log n-1/2)}{((n+1)\log n)}\}^{1/n}}
$$
Expanding ${\log \log(n+1)-1/2}$ into 
${\log \log n+\log[1+\frac{(\log(1+1/n))}{\log n}]}\\ {-1/2}$,
\begin{eqnarray}
{\{1+\frac{\log[1+\frac{\log(1+1/n)}{\log n}]}{(\log \log n-1/2)}\}<(n+1)^{1/n}(\log n)^{1/n}\times} \nonumber\\
 {\{1+\frac{(\log \log n-1/2)}{((n+1)\log n)}\}^{1/n}}, \label{eq3.10}
\end{eqnarray}
which we can prove as follows.
Let ${X=\frac{(\log \log n-1/2)}{\log[1+\frac{\log(1+1/n)}{\log n}]}>402}$ for ${n \geq 89}$, then the ${X}$th power of (\ref{eq3.10}) is 
\begin{equation} \label{eq3.11}
{\{1+\frac{1}{X}\}^{X}<(n+1)^{\frac{X}{n}}(\log n)^{\frac{X}{n}}\{1+\frac{(\log \log n-1/2)}{((n+1)\log n)}\}^{\frac{X}{n}}}.
\end{equation}

Trivially, ${\{1+\frac{1}{X}\}^{X}<e}$ for all ${X>0}$. However, since ${\log n>1}$ for ${n\geq 89}$, then ${(\log n)^{\frac{X}{n}}>1}$. The inequality ${\{1+\frac{(\log \log n-3/2)}{((n+1)\log n)}\}^{\frac{X}{n}}>1}$ is also true. Thus, it is sufficient to show that ${(n+1)^{\frac{X}{n}}>e}$ for all ${n\geq 89}$.

Since the denominator of the fraction ${\frac{X}{n}}$ is ${n\log[1+\frac{\log(1+1/n)}{\log n}]<\frac{1}{\log n}<1}$ for ${n\geq 3}$ (Proved in Section~\ref{sc2}),then
$${(n+1)^{\frac{X}{n}}>(n+1)^{(\log \log n-1/2)}>89}$$
for ${n\geq 89}$. Therefore,
$$
{\{1+\frac{1}{X}\}^{X}<e<89<(n+1)^{\frac{X}{n}}<(n+1)^{\frac{X}{n}}(\log n)^{\frac{X}{n}}\{1+\frac{(\log \log n-1/2)}{((n+1)\log n)}\}^{\frac{X}{n}}}.$$
Thus, the left side of (\ref{eq3.11}) is less than the right for ${n\geq 89}$, and (\ref{eq3.7}) and consequently (\ref{eq3.3}) and (\ref{eq3.2}) are proven. 

Therefore, we have proven (\ref{eq3.1}) and consequently (\ref{eq2.4}) for ${n\geq 89}$.
\end{enumerate}

\subsection{Proof of inequality (\ref{eq2.5})}

We can prove (\ref{eq2.5}) is true for all ${n\geq 195340}$ using the same method as for proving (\ref{eq2.4}). This is because the inequality ${(\log \log n-3/2)>1}$ holds, and so ${n\geq 195340>e^{e^{5/2}}}$

\subsection{Proof of inequality (\ref{eq2.6})}

To prove (\ref{eq2.6}) it is sufficient to show that 
\begin{equation} \label{eq3.12}
{\{m\log m+m\log \log m-m/2\}^{1/m}<\{n\log n+n\log \log n-3n/2\}^{1/n}}
\end{equation}
Holds for ${m\geq 3n}$ and ${n\geq 195340}$.

If we assume that ${m=3n}$, then we should show that
\begin{equation} \label{eq3.13}
{\{3n\log 3n+3n\log \log 3n-\frac{3n}{2}\}^{1/3n}<\{n\log n+n\log \log n-\frac{3n}{2}\}^{1/n}}
\end{equation}
for $n\geq 195340$.

The strictly decreasing properties of (\ref{eq2.4}) and (\ref{eq2.5})mean that (\ref{eq3.12}) is satisfied for all ${m>3n}$. Therefore, we need only show (\ref{eq3.12}) for ${m=3n}$. If we take the ${3n}$th power for (\ref{eq3.13}), then 
\begin{equation} \label{eq3.14}
{3n\log 3n+3n\log \log 3n-\frac{3n}{2}<\{n\log n+n\log \log n-\frac{3n}{2}\}^{3}}
\end{equation}
for $n\geq 195340$. Expanding the right side, 
$$
{(n\log n)^{3}+3(n\log n)^{2}\{n\log \log n-\frac{3n}{2}\}+3(n\log n)\{n\log \log n-\frac{3n}{2}\}^{2}}\\
+ {\{n\log \log n-\frac{3n}{2}\}^{3}}.
$$
Thus, we need to prove 
$$
{3n\log 3n+3n\log \log 3n-3n/2<(n\log n)^{3}+3(n\log n)^{2}\{n\log \log n-3n/2\}+}$$
$${3(n\log n)\{n\log \log n-3n/2\}^{2}+\{n\log \log n-3n/2\}^{3}}
$$
for ${n\geq 195340}$.

Note that ${3n\log 3n<(n\log n)^{3}}$, because if ${3n\log 3n=3n\log3+3n\log n<(n\log n)^{3}}$, then dividing each side by ${n(\log n)^{3}}$, 
$$
{n^{2}>\frac{3\log3}{(\log195340)^{3}}+\frac{3}{(\log195340)^{2}}>\frac{3\log3}{(\log n)^{3}}+\frac{3}{(\log n)^{2}}}
$$ 
for ${n\geq 195340}$, which holds. Thus, ${n^{2}>\frac{3\log3}{(\log n)^{3}}+\frac{3}{(\log n)^{2}}}$ for ${n\geq 195340}$, and ${3n\log 3n<(n\log n)^{3}}$ holds.

The second term of the inequality, 
$${3n\log \log 3n<3n(n\log n)^{2}\{\log \log n-3/2\}},$$is also satisfied.
The function ${z=(x\log x)^{2}-\log \log 3x}$ is strictly increasing. If ${z'=2x(\log x)^{2}+2x\log x-\frac{1}{(x\log 3x)}}$, then ${z'(195340)>0}$ and ${z'}$ is also strictly increasing due to ${x\log x}$, ${x(\log x)^{2}}$, and ${\frac{-1}{(x\log 3x)}}$ being strictly increasing for ${x\geq 195340}$. Therefore, ${z(x)>z(195340)>0}$ for all ${x>195340}$. 

Thus, ${(n\log n)^2>\log \log 3n}$ for ${n\geq 195340}$ and ${\frac{(\log n)^2}{\log \log 3n}>\frac{1}{n^{2}}}$. Since $\log \log n-3/2>1$ for $n\geq 195340$, then  ${\frac{(\log n)^{2}(\log \log n-3/2)}{\log \log 3n}>1/n^{2}}$ and ${(n\log n)^2(\log \log n-3/2)>}$ ${\log \log 3n}$. 

Multiplying both sides by ${3n}$, ${3n(n\log n)^2(\log \log n-3/2)>3n\log \log 3n}$ for ${n\geq 195340}$. Therefore, 
\begin{equation} \label{eq3.15}
{3n\log 3n+3n\log \log 3n<(n\log n)^{3}+3(n\log n)^{2}\{n\log \log n-3n/2\}}\end{equation}
for $n\geq 195340$. However, if we add ${-3n/2}$ to the left side and then add the two positive values ${3(n\log n)\{n\log \log n-3n/2\}^{2}}$ and ${\{n\log \log n-3n/2\}^{3}}$ for ${n\geq 195340}$ to the right side of (\ref{eq3.15}), then the inequality is satisfied and (\ref{eq3.14}), (\ref{eq3.13}), and (\ref{eq3.12}) are also satisfied. 

Therefore, we have proven the inequalities (\ref{eq2.4}) for ${n\geq 89}$ and (\ref{eq2.5}) and (\ref{eq2.6}) for ${n\geq 195340}$. 

From Silva \cite{T}, we know Theorem~1 holds for all ${n<195340}$, so that ${\{p_{n}\}_{n\leq 195340}<2770409<4 \times 10^{18}}$. Theorem~1 also says that ${m}$ is the smallest value greater than ${n}$ for every ${n\geq 1}$. The validity of this proposition may be concluded from (\ref{eq2.4}) and (\ref{eq2.5}), because these show that by progressing the value of ${n}$ to ${n+1}$, the inequalities' sides and the value of ${p_{n}^{1/n}}$ are decreasing to the first minimum of ${p_{m}^{1/m}}$, immediately before ${p_{n}^{1/n}}$ for ${m>n}$. Thus, ${m}$ should be the smallest number after ${n}$ so that ${p_{m}^{1/m}<p_{n}^{1/n}}$. Therefore, the first section of Theorem~1 is proved for ${n\geq 195340}$. Conversely, if ${p_{m}^{1/m}<p_{n}^{1/n}}$, then according to same Theorem~1, there is the smallest number, ${m_{1}}$, such that ${m_{1}>n}$ and ${p_{m_{1}}^{1/m_{1}}<p_{n}^{1/n}}$. Because, ${m_{1}}$, is the smallest value,thus ${m\geq m_{1}}$. If ${m=m_{1}}$, we have ${m=m_{1}>n}$ and the the conclusion stands. If ${m>m_{1}}$, then  ${p_{m}^{1/m}<p_{m_{1}}^{1/m_{1}}<p_{n}^{1/n}}$ and according to same Theorem~1, there is the smallest number ${m_{2}}$  so that ${p_{m}^{1/m}<p_{m_{2}}^{1/m_{2}}<p_{m_{1}}^{1/m_{1}}}$. This method can continue to ${k}$th value and ${m=m_{k}}$. Therefore, ${m\geq m_{k}>...>m_{1}>n}$ and finally, we have ${m>n}$, which completes the proof of Theorem~1. 

\setcounter{equation}{0}
\section{Theorem 2} \label{sc4}

\textit{Let ${\{p_{n}^{1/n}\}_{n\geq 1}}$ be an infinite sequence of real numbers, then this sequence is strictly decreasing }.

\textbf{\textit{Proof.}}
Using Theorem~1 and Lemma~3, we are able to prove Theorem~2. 
As shown in Theorem~1, we can prove that there is a natural number, ${m}$, for each natural number, ${n\geq 1}$, such that ${p_{m}^{1/m}< p_{n}^{1/n}}$, and ${m}$ is the first number greater than ${n}$ or ${m}$ is the smallest value after ${n}$. For Theorem~2, we need to prove that for each ${n\geq 1}$, there is an ${m}$ such that ${m=n+1}$. Silva and others have showed that Theorem~2 holds for all ${p_{n}}$'s with ${n\leq 195340 (p_{n}<2770409< 4\times 10^{18})}$. We need only show that it also holds for natural numbers ${n>195340}$.

Let the following sequence be true:
\begin{equation}\label{eq4.1}
{p_{195340}^{1/195340}<...<p_{2}^{1/2}<p_{1}^{1}}. 
\end{equation}  

From Theorem~1 and Lemma~3, there is the smallest value, ${m}$, such that ${m>195340}$ and 
\begin{equation}\label{eq4.2}
{p_{m}^{1/m}<p_{195340}^{1/195340}<...<p_{2}^{1/2}<p_{1}^{1}}.
\end{equation}

Consider the real number ${p_{195341}^{1/195341}}$, then regarding (\ref{eq4.2}), we have three cases:
\begin{enumerate}
\item The value of ${p_{195341}^{1/195341}}$ falls within the sequence defined in (\ref{eq4.1}). 

This case is impossible, because it contradicts the conditions of Theorem~1. If ${p_{195341}^{1/195341}}$ falls within (\ref{eq4.1}), then we may assume it arbitrarily falls between two values of this sequence as ${p_{195001}^{1/195001}}$ and ${p_{195000}^{1/195000}}$, so that $p_{195001}^{1/195001}<p_{195341}^{1/195341}<{p_{195000}^{1/195000}}$. Thus, from this inequality and Theorem~1, 195341 is the smallest number after 195000, such that ${p_{195341}^{1/195341}<p_{195000}^{1/195000}}$, but from the sequence (\ref{eq4.1}), this contradicts that the smallest value is 195001 and ${p_{195001}^{1/195001}<p_{195000}^{1/195000}}$. Therefore, ${p_{195341}^{1/195341}}$ cannot fall within (\ref{eq4.1}).

\item The value of ${p_{195341}^{1/195341}}$ exactly falls before ${p_{m}^{1/m}}$ such that
\begin{equation}\label{eq4.3}
{p_{195341}^{1/195341}<p_{m}^{1/m}<p_{195340}^{1/195340}<...<p_{2}^{1/2}<p_{1}^{1}}. 
\end{equation}

This case is also impossible. From Theorem~1, there is the smallest value, ${m>195340}$, such that ${p_{m}^{1/m}<p_{195340}^{1/195340}}$, from (4.3), we have ${p_{195341}^{1/195341}<p_{m}^{1/m}<p_{195340}^{1/195340}}$ 
 and from this inequality, Theorem~1, and Lemma~3, we have, ${195340<m<195341}$, which is impossible, thus the smallest value before ${p_{195340}^{1/195340}}$ is the value ${p_{195341}^{1/195341}}$, not ${p_{m}^{1/m}}$.

\item The value ${p_{195341}^{1/195341}}$ exactly falls between ${p_{m}^{1/m}}$ and ${p_{195340}^{1/195340}}$, 
\begin{equation}\label{eq4.4}
{p_{m}^{1/m}<p_{195341}^{1/195341}<p_{195340}^{1/195340}<...<p_{2}^{1/2}<p_{1}^{1}}.
\end{equation}
This is the only case which can occur.
\end{enumerate}

Arguing similarly for the sequence ${\{p_{n}^{1/n}\}_{n>195341}}$ completes the proof of Theorem~2, and proves Firoozbakht's conjecture. 

\textbf{Acknowledgment}

The author would like to thank Prof. Geoffrey Campbell for endorsing me for submitting the paper to arXiv.org. 


\begin{thebibliography}{99}
\bibitem{Fe} J. Feliksiak,The symphony of primes,distribution of primes and Riemann's hypothesis, Xlibris,(2013) 34-42
\bibitem{C} H. Cramer,On the order of magnitude of the difference between consecutive prime numbers,Acta Arith.\textbf{2}(1936) 23-46
\bibitem{S} D. Shanks ,On Maximal Gaps between Successive Primes, Math.Comput.\textbf{18},88 (1964) 646-651
\bibitem{N} T.R. Nicely,New maximal prime gaps and first occurrences, Math.Comput.\textbf{68},227(1999) 1311-1315
\bibitem{F} C. Rivera,Conjecture 30. The Firoozbakht Conjecture, www.primepuzzles.net (2012)
\bibitem{T} T. Oliviera e Silva, Gaps between consecutive primes, sweet.ua.pt/tos/gaps.html (2015)
\bibitem{K} A. Kourbatov, Upper Bounds for Prime Gaps
Related to Firoozbakht’s Conjecture, J. Integer Sequence\textbf{18}(2015) 1-7
\bibitem{RS}J.B. Rosser, L.Schoenfeld,Approximate formulas for some functions of prime numbers, Illinois J.Math.\textbf{6} (1962) 64-94
\end{thebibliography}
\end{document}